\documentclass[12pt,reqno]{amsart}

\usepackage{amsmath}
\usepackage{amscd}
\usepackage{amssymb}
\usepackage{latexsym}
\usepackage{MnSymbol}
\usepackage{undertilde}
\usepackage{amsrefs}
\usepackage{nicefrac}

\input xy
\xyoption{all} 
\newdir{ >}{{}*!/-5pt/@{>}}

\theoremstyle{plain}
\newtheorem{theorem}{Theorem}[section]
\newtheorem*{theorem*}{Theorem}
\newtheorem{lemma}[theorem]{Lemma}
\newtheorem*{lemma*}{Lemma}
\newtheorem{corollary}[theorem]{Corollary}
\newtheorem{proposition}[theorem]{Proposition}
\newtheorem{defprop}[theorem]{Proposition-Definition}

\theoremstyle{remark}
\newtheorem{remark}[theorem]{Remark}

\theoremstyle{definition}
\newtheorem{definition}[theorem]{Definition}
\newtheorem{example}[theorem]{Example}

\makeatletter
\def\revddots{\mathinner{\mkern1mu\raise\p@
\vbox{\kern7\p@\hbox{.}}\mkern2mu
\raise4\p@\hbox{.}\mkern2mu\raise7\p@\hbox{.}\mkern1mu}}
\makeatother 
\newcommand{\bgl}{\begin{equation}} 
\newcommand{\egl}{\end{equation}}
\newcommand{\bgloz}{\begin{equation*}} 
\newcommand{\egloz}{\end{equation*}}
\newcommand{\bgln}{\begin{eqnarray}} 
\newcommand{\egln}{\end{eqnarray}}
\newcommand{\bglnoz}{\begin{eqnarray*}} 
\newcommand{\eglnoz}{\end{eqnarray*}}
\newcommand{\btheo}{\begin{theorem}}
\newcommand{\etheo}{\end{theorem}}
\newcommand{\btheooz}{\begin{theorem*}}
\newcommand{\etheooz}{\end{theorem*}}
\newcommand{\blemma}{\begin{lemma}}
\newcommand{\elemma}{\end{lemma}}
\newcommand{\blemmaoz}{\begin{lemma*}}
\newcommand{\elemmaoz}{\end{lemma*}}
\newcommand{\bproof}{\begin{proof}}
\newcommand{\eproof}{\end{proof}}
\newcommand{\bbew}{\begin{beweis}}
\newcommand{\ebew}{\end{beweis}}
\newcommand{\bremark}{\begin{remark}}
\newcommand{\eremark}{\end{remark}}
\newcommand{\bex}{\begin{example}\em}
\newcommand{\eex}{\end{example}}
\newcommand{\bdefin}{\begin{definition}}
\newcommand{\edefin}{\end{definition}}
\newcommand{\bprop}{\begin{proposition}}
\newcommand{\eprop}{\end{proposition}}
\newcommand{\bdefprop}{\begin{defprop}}
\newcommand{\edefprop}{\end{defprop}}
\newcommand{\bcor}{\begin{corollary}}
\newcommand{\ecor}{\end{corollary}}
\newcommand{\bfa}{\begin{cases}} 
\newcommand{\efa}{\end{cases}}
\newcommand{\n}{\par\noindent}

\newcommand{\mn}{\par\medskip\noindent}

\newcommand{\cK}{\mathcal K}
\newcommand{\cL}{\mathcal L}

\def\Cz{\mathbb{C}}

\def\Nz{\mathbb{N}}

\def\Tz{\mathbb{T}}
\def\Zz{\mathbb{Z}}
\def\Z{\mathbb{Z}}

\def\1z{\mathbb{1}}
\newcommand{\lori}{\longrightarrow}

\def\SEMI{\mbox{$\times\kern-2pt\vrule height5pt width.6pt \kern3pt $}}

\newcommand{\comment}[1]{}  

\newcommand*{\id}{\mathrm{id}}

\begin{document}

\title[]{Semigroup C*-algebras and toric varieties}
\author{Joachim Cuntz}
\email{cuntz@uni-muenster.de}
\address{Mathematisches Institut, Einsteinstr. 62, 48149 M\"unster, Germany}
\date{\today}

\begin{abstract}
Let $S$ be a finitely generated subsemigroup of $\Zz^2$. We derive a general formula for the $K$-theory of the left regular C*-algebra $C^*_\lambda S$.
\end{abstract}
\maketitle
\section{Introduction}
Let $S$ be a finitely generated subsemigroup of $\Zz^n$. Then its monoid algebra $\Cz S$ is a finitely generated $\Cz$-algebra with no non-zero nilpotent elements. It is therefore the coordinate ring of an affine variety over $\Cz$. Such varieties are called affine toric varieties (they carry an action of an $n$-dimensional torus). Of course, here we may replace $\Cz$ by an arbitrary field. General references for toric varieties and the corresponding semigroups are for instance \cite{CLS} or \cite{Neeb}.

In this article we study the left regular semigroup C*-algebra $C^*_\lambda S$ which, in contrast to $\Cz S$, is generated not only by the elements of $S$ but also by their adjoints (it is in fact generated by the enveloping inverse semigroup to $S$). It also carries a natural action of $\Tz^n$.

As we shall see (Lemma \ref{L10}), the case $n=1$ is without interest: for a non-trivial subsemigroup $S$ of $\Zz$ the C*-algebra $C^*_\lambda S$ is in fact always isomorphic (non-canonically) to the ordinary Toeplitz algebra $C^*_\lambda \Nz$. But our main objective in the present paper is the computation of the $K$-theory of $C^*_\lambda S$ for a finitely generated subsemigroup $S$ of $\Zz^2$. In \cite{CEL1}, \cite{CEL2} we had determined the $K$-theory of a large class of semigroup C*-algebras using the independence condition introduced by Xin Li. The interesting feature of the semigroups that we meet here, however is that they do not satisfy this condition except in trivial cases. One consequence is that the $K$-theory contains a torsion part and, as another consequence, we cannot rely here on an elegant general method for its determination.

But, by a more detailed study of the structure of arbitrary finitely generated subsemigroups of $\Zz^2$, we are able to show that the $K$-theory of the C*-algebra is always described by a simple formula involving only the `faces' of the semigroup.

In \cite{PaSc}, on the basis of previous results in \cite{JiKa}, \cite{JiXia} and \cite{Park}, a formula for the $K$-theory of $C^*_\lambda S$ (which looks different from ours, but gives the same result) had been established in the important special case of a `saturated' finitely generated subsemigroup of $\Zz^2$. These authors consider `Toeplitz algebras' associated with cones, but, upon inspection, their Toeplitz algebra is exactly the left regular C*-algebra of the (automatically saturated) semigroup defined by a cone in $\Zz^2$. In the saturated case our computation is somewhat more direct than the one in \cite{PaSc}. But much of the analysis in our paper is really concerned with the non-saturated case. The result in \cite{PaSc} actually also covers saturated subsemigroups of $\Zz^2$ which are defined by a cone where one face has irrational slope. Such semigroups are not finitely generated. But the semigroup defined by a half-plane with irrational slope in $\Zz^2$ does satisfy the independence condition. Thus one easily derives the main result in \cite{JiKa}, on the $K$-theory of the corresponding Toeplitz algebra, by applying the general method based on independence from \cite{CEL1}. Using this, the case of a cone with one irrational face considered in \cite{PaSc} also lies within the scope of our methods.

Even though the exact sequences in $K$-theory, that we use, become significantly more complicated in higher dimensions, it may well be possible that our argument can be extended to subsemigroups of $\Zz^n$ for $n>2$.
\section{Toric varieties}
Much of the material in this section is well known. We consider a finitely generated subsemigroup $S$ of $\Zz^n$. We write the semigroup operation as addition and always assume that a semigroup contains $0$. It is then easily seen that the subgroup $S-S$ generated by $S$ in $\Zz^n$ is the enveloping group of $S$ and is of course isomorphic to $\Zz^k$ for some $k$. Thus, without restriction of generality, we may assume that $\Zz^n$ already is the enveloping group and therefore that $S$ generates $\Zz^n$ and that this embedding is natural. We will also assume that $S\cap (-S) = \{0\}$, i.e that $S$ contains no invertible elements besides $0$.

Summarizing, we assume from now on that $S$ is a finitely generated subsemigroup of $\Zz^n$, for some $n$, generating $\Zz^n$ as a group. We also assume that $S\cap (-S) = \{0\}$. The following easy lemma will be used in two places
\blemma\label{Le} Let $Y=\{y_1,\ldots ,y_m\}$ be a finite set in $\Zz^n$. Then there is $z\in S$ such that $z+Y \subset S$.
\elemma
\bproof
Let $x_1,\ldots ,x_l$ denote the generators of $S$.
Since $S$ generates $\Zz^n$ as a group, there are $k^i_j\in \Zz$ such that $y_i=\sum_j k^i_jx_j$. We may assume that the $k^i_j$ are ordered in such a way that $k^i_j< 0$ for $i\leq r_i$ and $k^i_j\geq 0$ for $i> r_i$. Denote by $\bar{y}_i = \sum_{j=1}^{r_i}k^i_jx_j$ the `negative part' of $y_i$. Then $z=-\sum_i \bar{y}_i$ has the required property.
\eproof

A subsemigroup $F\subset S$ is said to be a face of $S$ if $x+y \in F$ with $x,y\in S$ implies that $x,y \in F$. A semigroup $S\subset \Zz^n$ is said to be saturated if $kx\in S,\, x\in\Zz^n,\, k\in \Nz\setminus \{0\}$ implies that $x\in S$. For an arbitrary subsemigroup $S$ of $\Zz^n$ let $\bar{S}$ denote the saturation of $S$, i.e. the semigroup consisting of all $s\in \Zz^n$ for which an integral multiple $ks$, $k\in \Nz\setminus \{0\}$ lies in $S$. There is a bijection between faces in $S$ and faces in $\bar{S}$ as follows
\blemma\label{L9}(see also \cite{Neeb} Lemma II.7) Let $F\subset S$ be a face of $S$. Then $\bar{F}$ is a face in $\bar{S}$ and $F=\bar{F}\cap S$. Conversely, if $G$ is a face in $\bar{S}$, then $G\cap S$ is a face in $S$.\elemma
\bproof Assume that $F$ is a face and let $x',y'$ be in $\bar{S}$ and $f'$ in $\bar{F}$ such that $x'+y' = f'$. We have that $kx',jy',nf'$ are in $S$ and $F$ respectively for suitable $k,j,n$ in $\Nz$. Let $m$ be the least common multiple of $k,j,n$. It follows that
$mx'+my' = mf'$ with $mx', my' \in S$ and $mf'\in F$. Since $F$ is a face, $mx'$ and $my'$ are in $F$ and thus $x',y'$ are in $\bar{F}$. By the defining property of a face we have that $F = \bar{F}\cap S$.

Conversely, assume that $\bar{F}$ is a face. If $x,y\in S$ and $f\in \bar{F}\cap S$ are such that $x+y=f$, then $x,y$ are in $\bar{F}$ and thus also in $F=\bar{F}\cap S$. \eproof
A subsemigroup $T$ of $S$ is said to be one-dimensional, if the subgroup of $\Zz^2$ it generates, is isomorphic to $\Zz$.
Since a saturated generating semigroup of $\Zz^2$ is determined by a convex cone (see e.g. \cite{Neeb} Lemma II.7), Lemma \ref{L9} implies that any generating subsemigroup $S$ of $\Zz^2$ has exactly two one-dimensional faces. Note, that a subsemigroup of $\Zz$ that contains no non-zero invertibles has to lie entirely in $\Nz$ or in $-\Nz$. Therefore the structure of one-dimensional subsemigroups of $S$ is determined by the following Lemma.
\blemma\label{L8} Let $F$ be a finitely generated subsemigroup of $\Nz$. Then there is $d\in \Nz$ such that $F\subset d\Nz$ and such that $d\Nz \setminus F$ is finite.\elemma
\bproof Let $F'= F-F$ be the subgroup of $\Zz$ generated by $F$. Then there is $d\in \Nz$ such that $F'=d\Zz$. It follows that $F\subset d\Nz$. Let $m\in \Nz$ such that $md\in F$. By Lemma \ref{Le} there is $z\in F$ such that $z+jd$ is in $F$ for $j=1\ldots ,m-1$. Since $md\in F$, it follows that $z+\Nz d$ is contained in $F$.
\eproof
Let $F \subset S$ be a one-dimensional subsemigroup. By Lemma \ref{L8} there exists a unique $a \in \Zz^n$ such that $F$ is contained in $\Nz a$ with finite complement. Moreover $F$ then generates $\Zz a$ as a group.
\bdefin\label{ag}  Given $F\subset S$ as above, we say that the element $a$ is the asymptotic generator of $F$.\edefin
Recall that the quotient of a commutative semigroup $S$ by a subsemigroup $F$ is the semigroup consisting of equivalence classes of elements $s$ in $S$ for the equivalence relation $s_1\sim s_2 \;\Longleftrightarrow\; \exists f_1,f_2\in F$ such that $s_1+f_1=s_2+f_2$.
\blemma\label{L1} Let $F \subset S$ be as above and $a$ the asymptotic generator of $F$.
Denote by $x\mapsto \dot{x}$ the quotient map $S\to S/F$.
Then $\dot{x}=\dot{y}$ for $x,y\in S$ if and only if $(x+\Zz a)\cap S = (y+\Zz a)\cap S$. If $\dot{x}\neq\dot{y}$, then $(x+\Zz a)\cap (y+\Zz a)= \emptyset$.\elemma
\bproof If $\dot{x}=\dot{y}$, then there are $f_1,f_2$ in $F$ such that $x+f_1=y+f_2$ and thus that $x-y =f_2 -f_1 \in\Zz a$. This implies that $x\in y+\Zz a$ and $y\in x+\Zz a$.

Conversely, assume that $x=y+ka$ with $k\in \Zz$. By Lemma \ref{L8} there is $n\in \Nz$ such that $ka+na$ and $na$ are in $F$. It follows that $x+na = y+ka+na$ and thus that $\dot{x}=\dot{y}$.

The same argument shows that, if $x+k_1a=y+k_2a$ for $k_1,k_2\in \Zz$, then $\dot{x}=\dot{y}$.\eproof

\bcor\label{C1} Let $S$ and $F$ be as in Lemma \ref{L1}. Then $S$ is a disjoint union $$S =\bigsqcup_{\dot{x}\in S/F} (x+\Zz a)\cap S $$\ecor
\bproof Since, for $\dot{x}\in S/F$, the set $x+ \Zz a$ does not depend on the representative $x$, this is an immediate consequence of Lemma \ref{L1}.\eproof

\blemma Let $F$ be a non-trivial face in $S$. Then for each $x\in S$ with $x\notin F$, we have that $F\subset S \setminus (S+x)$. \elemma
\bproof If $F$ is a face, then $x\notin F$ implies that $F\cap (S+x) = \emptyset$.
\eproof

\blemma\label{L4} Let $x\in S$ and $\langle x \rangle$ the subsemigroup generated by $x$. Then the quotient map $y \mapsto \dot{y} \in S/\langle x \rangle$ induces a bijection between $S\setminus (S +x)$ and $S/\langle x \rangle$. \elemma
\bproof Assume that $\dot{y}_1 = \dot{y}_2$ for $y_1, y_2 \in S\setminus (S+x)$. Then there is $n\in \Nz$ such that $y_1 +nx =y_2$ or $y_2 +nx =y_1$. Since $y_1, y_2 \in S\setminus (S+x)$, $n$ has to be zero, so that $y_1 = y_2$. This shows injectivity.

To show surjectivity, take $\dot{y} \in S/\langle x \rangle$ represented by $y\in S$. There is a minimal $n\in \Nz$  such $y\in S+nx$, i.e. $y=z+nx$ and $z\notin S+x$. Then $\dot{z} =\dot{y}$.
\eproof
\blemma\label{L5} Let $S\subset \Zz^2$ and let $F_1,F_2$ denote the two one-dimensional faces of $S$. Let $a_1$, $a_2$ be the asymptotic generators of $F_1$, $F_2$ and $C$ the cone in $\Zz^2$ spanned by $a_1$ and $a_2$ (i.e. $C=\bar{S}$ using the notation above). Then there is $z\in S$ such that $z + C\subset S$.
\elemma
\bproof Let $x_1,\ldots ,x_n$ denote the generators of $S$ and let $b_1,b_2$ be multiples of $a_1,a_2$ such that $b_i\in F_i$, $i=1,2$. Let $P=\{y_1,\ldots ,y_m\}$ denote the set of all elements in $\Zz^2$ that lie inside the parallelogram spanned by $b_1$ and $b_2$. Then $P+F'=C$ for the subsemigroup $F'$ of $S$ spanned by $b_1,b_2$.
By Lemma \ref{Le} there is $z\in S$ such that $z+P\subset S$. Then, since $C=F'+P$, also $z+C\subset S$.
\eproof
\blemma\label{L7} Let $F\subset S$ be a one-dimensional face and $a$ the asymptotic generator of $F$. Then $F= S\cap \Zz a$.\\
For each $x\in S$, $S\setminus (S+x)$ is a finite union of finitely many translates of $F_1$ and $F_2$, and of a finite set.
\elemma
\bproof It is clear that $F\subset S\cap \Zz a$. Conversely, let $ka \in S$ for $k\in \Zz$. By Lemma \ref{L8} there is $n$ in $\Nz$ such that $na$ and $na +ka$ are in $F$. Since $F$ is a face, this implies that $ka\in F$.

Let $C$ be as in \ref{L5}. If $x\in S$, then by Lemma \ref{L5}, there is $z\in S$ such that $C+z\subset S+x$. Now, $S\setminus (C+z)$ is a finite union of subsets of the form $(y+\Zz a_1)\cap S$ or $(y+\Zz a_2)\cap S$ (each diagonal in $\Zz^2$ parallel to $a_i$ is a finite union of subsets of the form $y+\Zz a_i$) - and thus, up to a finite set, a finite union of translates $y+F_i$, $i=1,2$. Therefore also $S\setminus (S+x)$ is a finite union of subsets $(y+F_i)\cap (S\setminus (S+x))$, $i=1,2,\, y\in S$. \\
By Corollary \ref{C1}, for each translate $y+F_i$, the intersection with $S+x$ is empty or has finite complement in $y+F_i$.
\eproof

\section{The regular C*-algebra for a toric semigroup}
We consider a finitely generated generating subsemigroup $S$ of $\Zz^2$ and denote by $F_1$, $F_2$ the two one-dimensional faces of $S$. We denote by $\lambda$ the left regular representation of $S$ on $\ell^2 S$ and by $C^*_\lambda S$ the C*-algebra generated by $\lambda (S)$.
As usual, there is the commutative sub-C*-algebra $D$ of $C^*_\lambda S$ which is generated by all range projections of the partial isometries obtained as all possible products of the $\lambda (s), s\in S$ and their adjoints.
\blemma $D$ contains all orthogonal projections onto $\ell^2(X)$ where $X$ is a finite subset of $S$. Consequently, $C^*_\lambda S$ contains the algebra $\cK$ of all compact operators on $\ell^2S$.
\elemma
\bproof $D$ contains the orthogonal projection onto $\ell^2(X)$ where $X=(S\setminus (S+f_1))\cap (S\setminus (S+f_2))$, $f_1\in F_1, f_2\in F_2$. Lemma \ref{L7} implies that $X$ is finite. Consider now all subsets of $S$ obtained as the intersection of $X$ with finitely many translates $s+X$ with $s\in S-S$. Let $Y$ denote a minimal set in this family. Then $(Y+s_1)\cap (Y+s_2)=\emptyset$ whenever $s_1\neq s_2$. Let now $y_1,y_2\in Y$. Then, since the enveloping group $S-S$ is $\Zz^2$, there are $s_1,s_2\in S$ such that $y_1+s_1=y_2+s_2$. This implies $s_1=s_2$ and thus also $y_1=y_2$. We see that $Y$ consists of only one point. The one-dimensional projection onto $\ell^2 (Y)$ is in $D$ and therefore also all of its translates.
\eproof
Given a subset $X$ of $S$, denote by $e_X$ the orthogonal projection onto the subspace $\ell^2X\subset \ell^2 S$.
\blemma $D$ is generated by the projections of the form $\lambda(s)e_F\lambda (s)^*$, for $s\in S$ and $F$ a face of $S$. \elemma
\bproof We show first that $e_F$ is in $D$ for each face. This is clear for the trivial faces $\{0\}$ and $S$. Thus let $F$ be one of the two one-dimensional faces and $C$ and $z$ as in Lemma \ref{L5}. There is $d\in \Zz^2$ such that $C\setminus ((C+d)\cap C)$ is the face of $C$ which contains $F$. Replacing $z$ by a translate $z+x$, for a suitable $x$ if necessary, we may clearly assume that $z+d+C\subset S$. It follows that $F= ((z+S))\setminus(z+d+S))-z) \cap S$ and thus that $$e_F= \lambda (z)^*\left(\lambda (z)\lambda (z)^*-\lambda (z+d)\lambda (z+d)^*\right)\lambda (z)$$

Denote by $D_0$ the subalgebra generated by all projections $\lambda(s)e_F\lambda (s)^*$. Since $e_F$ is in $D$, we have $D_0\subset D$. Moreover $D_0$ then contains all diagonal projections of finite rank. Lemma \ref{L5} also implies that the complement of any range projection of $\lambda (s)$ for $s\in S$ is a linear combination of finitely many translates of projections of the form $e_F$. This shows that $D_0 =D$. \eproof
\bremark Toric semigroups typically do not satisfy the independence condition which says that the projections in $D$, obtained as range projections of products of elements $\lambda (s),\, s\in S$ and their adjoints, should be linearly independent. As a simple example consider the semigroup $S\subset \Zz^2$ defined by the cone spanned by the vectors $(2,1)$ and $(2,-1)$. Then the intersection of $(2,1)+S$ and $(2,0)+S$ equals the union of $(4,1)+S$ and $(4,0)+S$. This kind of phenomenon occurs for all toric semigroups except for the trivial ones.\eremark
\blemma\label{Lq} Let $F$ be a two-dimensional subsemigroup of $S$. Then the quotient $S/F$ is a finite abelian group and equal to the quotient $(S-S)/(F-F)$ of the enveloping groups. If $a=(k,l)$ and $b=(m,n)$ are generators of $F-F$, then the number of elements in $S/F$ is given by the absolute value of the determinant $$\det\left(\begin{array}{cc}k&m\\l&n \end{array}\right)$$
\elemma
\bproof Elements $x,y$ in $F$ become equal in $F/S$ if and only if there are $f,g$ in $F$ such that $x+f=y+g$ and thus if and only if $x-y=g-f$, i.e. iff $x,y$ become equal in $(S-S)/(F-F)$. This means that the map $S/F \to (S-S)/(F-F)$ is injective. Now, $F-F$ is a two-dimensional subgroup of $S-S=\Zz^2$ and therefore $(S-S)/(F-F)$ is finite. Thus the image of $S/F$ in $(S-S)/(F-F)$ is a subsemigroup of a finite group and therefore already a group. The formula for the number of elements in $S/F$ is well known and follows from the elementary divisor theorem.
\eproof
\blemma\label{L10} Let $F$ be a finitely generated subsemigroup of $\Nz$ generating $\Nz$. Then $C^*_\lambda F\cong C^*_\lambda \Nz$. Moreover, viewed as subalgebras of $\cL (\ell^2\Nz )$,  the algebra $C^*_\lambda F$ is a subalgebra of $C^*_\lambda \Nz$ such that $C^*_\lambda F /\cK = C^*_\lambda \Nz/\cK$.
\elemma
\bproof Lemma \ref{L8} shows that $M=\Nz\setminus F$ is a finite set. Let $n\in \Nz$ be large enough so that $M\subset \{0,\cdots ,n\}$ and $e_n$ be the projection onto $\ell^2 \{0,\cdots ,n\}$. Then we can find $f,g\in F$ such that $\lambda_\Nz (1) (1-e_n) = \lambda_F (g)^*\lambda (f)(1-e_n)$ (we denote here by $\lambda_\Nz$, $\lambda_F$ the left regular representations on $\ell^2\Nz$ and $\ell^2F$, respectively). It follows that $C^*_\lambda F = (1-e_M)\,C^*_\lambda (\Nz)\,(1-e_M)$. Moreover, using the fact that $C^*_\lambda \Nz$ is the universal C*-algebra generated by a single isometry, it is trivially seen that $C^*_\lambda (\Nz)\cong (1-e_M)\,C^*_\lambda (\Nz)\,(1-e_M)$.
\eproof
Let $F_1$, $F_2$ be the two one-dimensional faces of $S$ and denote by $I_1,I_2$ the closed ideals generated in $C^*_\lambda S$ by $e_{F_1}$, and by $e_{F_2}$, respectively.
\blemma The intersection $I_1\cap I_2$ is equal to $\cK (\ell^2 S)$. Each quotient $I_j/\cK(\ell^2 S)$ is isomorphic to $\cK (\ell^2 (S/F_j))\otimes C(\Tz)$. \\ Moreover, the quotient  $C^*_\lambda S/(I_1 +I_2)$ is isomorphic to $C(\Tz ^2)$.
\elemma
\bproof The first assertion follows from the fact that each intersection of a translate of $F_1$ and a translate of $F_2$ contains at most one point. The second assertion is a consequence of Corollary \ref{C1}  in combination with Lemma \ref{L7} and Lemma \ref{L10}. Finally, Lemma \ref{L7} also shows that any element $\lambda(s)$, $s\in S$ becomes unitary in the quotient $C^*_\lambda S/(I_1 +I_2)$ so that the quotient is isomorphic to the C*-algebra of the enveloping group $\Zz^2$ of $S$.
\eproof
As customary, we will, from now on, not distinguish between the algebras of compact operators on different separable infinite-dimensional Hilbert spaces and just write $\cK$. For the $K$-theory of the C*-algebra $C^*_\lambda S/\cK$ we obtain the following six-term exact sequence
\bgl\label{PV}\xymatrix{K_*(\dot{I}_1)\oplus K_*(\dot{I}_2))\ar[r]&
K_*(C^*_\lambda S/\cK)\ar[r]&K_*(C\Tz^2)\ar@/^9mm/[ll]}\egl\mn
where $\dot{I}_j$ denotes the quotient $I_j/\cK$.
\blemma\label{L17} Let $a_1,a_2$ be the asymptotic generators of the faces $F_1,F_2$ ordered in such a way that $\det (a_1,a_2)$ is positive (this implies that $\det (a_1,s)$ and $\det (s,a_2)$ are positive for all $s\in S$). Denote by $\pi$ the quotient map $C^*_\lambda S\to C(\Tz^2)$. Let $a_j=(x_j,y_j)$ and let $s=(m,n)\in S\setminus F_j$. The index map $K_1(C\Tz^2) = K_1(C^*_\lambda S)/(I_1+I_2)) \to K_0(\dot{I}_j)\cong\Zz$, for the extension \eqref{PV}, maps the class of $\pi(\lambda (s))$ to
$$(-1)^{j+1}\det (a_j,s)= (-1)^{j+1}\det\left(\begin{array}{cc}x_j&m\\y_j&n \end{array}\right)$$
\elemma
\bproof By Lemma \ref{L7}, the set $S\setminus (S+s)$ which represents the index of $\pi (\lambda (s))$ is, up to finite subsets, a union of finitely many translates of $F_1$ and $F_2$. When we project to $\dot{I}_j$, the number of the translates of $F_j$, that we obtain, is, according to Lemma \ref{L4}, given by the number of elements in $S/(F_j+\langle s\rangle )$ which in turn by Lemma \ref{Lq} is determined by the absolute value of the determinant above.
\eproof
\blemma\label{L19} Let $C$ and $z$ be as in Lemma \ref{L5} and let $a_1,a_2$ denote the asymptotic generators of the faces $F_1,F_2$. Then the projection $E$ onto $\ell^2 (z+C)$ is in $C^*_\lambda S$. The formulas $v_i=\lambda(a_i)E$, $i=1,2$ define elements in $C^*_\lambda S$. The isometries $v_1,v_2$ are relatively prime in the sense that $$v_1v_2^* = v_2^*v_1$$ In particular, the C*-subalgebra of $C^*_\lambda S$ generated by $v_1, v_2$ is isomorphic to the Toeplitz algebra $C^*_\lambda\Nz^2$.\elemma
\bproof Each diagonal in $S$ parallel to $a_i$ is invariant under addition of $F_i$ and therefore, up to finite sets, a finite union of translates of $F_i$. It follows that the complement of $z+C$ in $S$ is, up to finite sets, a finite union of translates of $F_1$ and of $F_2$. Since the projection onto $\ell^2$ of such a translate is in $C^*_\lambda S$, we see that $E=e_{z+C}$ is in $C^*_\lambda S$.

Consider now $v_1$ and $v_2$. Using Lemmas \ref{L8} and \ref{L10} in combination with the fact that $z+C$ is invariant under addition of $a_i$, we may assume that the formulas for $v_1,v_2$ actually define elements of $C^*_\lambda S$ (if $x\in F_i$ such that $x+a_i \in F_i$, then $\lambda(a_i)E = \lambda (x)^*\lambda(x+a_i)E$). The primeness condition is equivalent to the fact that the range projection of the product $v_1v_2$ is equal to the product of the range projections of the $v_1$, $v_2$. Therefore we have to show that $(z+C+a_1)\cap(z+C+a_2)= (z+C+a_1+a_2)$.
But, since $a_1$, $a_2$ span the boundary of the cone, one clearly has that $(C+a_1)\cap (C+a_2) = (C+a_1+a_2)$.
\eproof
Let $S$ be a finitely generated subsemigroup of $\Zz^2$ generating $\Zz^2$ as a group. Let $a_1=(x_1,y_1)$ and $a_2=(x_2,y_2)$ denote the asymptotic generators of the two one-dimensional faces $F_1$ and $F_2$ of $S$. In the following we use the integral $2\times 2$ -matrices
$$M\,=\, \left(\begin{array}{cc}y_2&-x_2\\-y_1&x_1 \end{array}\right)\qquad M^\perp\,=\, \left(\begin{array}{cc}x_1&x_2\\y_1&y_2 \end{array}\right)$$
Here again we order $a_1,a_2$ so that $\det M^\perp$ is positive. Note that $M$ is the adjugate matrix to $M^\perp$ in the sense of Cramer's rule so that $\det M = \det M^\perp$ and $MM^\perp = \det M \,\textbf{1}$.
\blemma\label{L18} Consider the extensions $0\to I\to C^*_\lambda S \to C(\Tz^2)\to 0$ and $0\to I'\to C^*_\lambda \Nz^2 \to C(\Tz^2)\to 0$, where $I$, $I'$ denote the kernels of the quotient maps. By Lemma \ref{L19} there is a natural map $\kappa: C^*_\lambda \Nz^2 \to C^*_\lambda S$ which maps the generators of $C^*_\lambda \Nz^2 $ to $v_1$, $v_2$, where $v_1$, $v_2$ are as in Lemma \ref{L19}.
Then we have the following:
\begin{enumerate}
  \item $K_0(I)=K_0(I')=\Z^2$ and $K_1(I)=K_1(I')=\Zz$. The generator of $K_1(I)$ is represented by $w=\lambda (a_1)e_{F_1} + \lambda (a_2)^*e_{F_2}$ (this is unitary mod $\cK$).
  \item The map $K_0(C\Tz^2)\to K_0(C\Tz^2)$ induced by $\kappa$ maps the Bott element $b$ to $(\det M)\, b$.
  \item The boundary map $K_1(C\Tz^2)= \Zz^2 \to K_0(I)=\Zz^2$ is given by multiplication by $M$.
  \item The map $K_1I'\cong \Zz \to K_1(I)\cong \Zz$ induced by $\kappa$ is multiplication by $\det M$.
\end{enumerate}
\elemma
\bproof (1) It follows from Corollary \ref{C1} and Lemma \ref{Lq} that the ideals $I$ and $I'$ are stably isomorphic (and thus, since both are stable, even isomorphic). In the long exact $K$-theory sequence for the extension $0\to I'\to C^*_\lambda \Nz^2 \to C(\Tz^2)\to 0$ we know that $K_0(C^*_\lambda \Nz^2)=\Zz$ with generator $[1]$ and that $K_1(C^*_\lambda \Nz^2)=0$. This shows that $K_1(I) \cong \Zz$. The fact that the generator is represented by $w$ follows from the $K$-theory sequence for the extension $0\to \cK \to I \to I/\cK \to 0$ and the fact that $I/\cK \cong (\cK \otimes C(\Tz))\oplus (\cK \otimes C(\Tz))$.

(2) It is obvious by definition that $\kappa_* : K_1(C(\Tz^2))\cong \Zz^2 \lori K_1(C(\Tz^2))\cong \Zz^2$ is given by multiplication by the matrix $M^\perp$. The Bott element is represented by the exterior product of the generators of $K_1(C(\Tz^2))\cong \Zz^2$. The map induced by multiplication by $M$ on the exterior product is $\det M$ by definition of the determinant.

(3) In the isomorphism $K_0(\dot{I})\cong \Zz^2$ we identify $K_0(\dot{I}_2)$ with the first component and $K_0(\dot{I}_1)$ with the second component of $\Zz^2$ (this convention is used for the identification of the maps in diagram \eqref{mu} below). By Lemma \ref{L17} (and keeping in mind the reverse identification of the components in $\Zz^2$) we know that the boundary map for the extension \eqref{PV} maps an element $s=(m,n) \in K_1(C\Tz^2)\cong \Zz^2$ to the element $(k_2,k_1)$ in $\Zz^2\cong K_0(\dot{I})$ with components
$$k_i = (-1)^{i+1}\det\left(\begin{array}{cc}x_i&m\\y_i&n \end{array}\right)$$
where $(x_i,y_i)$ are the components of the asymptotic generators $a_i$ of $F_i$, $i=1,2$. In other words
$$\left(\begin{array}{c}k_2\\k_1 \end{array}\right)\; =\; \left(\begin{array}{cc}y_2&-x_2\\-y_1&x_1 \end{array}\right)\,\left(\begin{array}{c}m\\n \end{array}\right)$$
This describes the boundary map to $K_0(\dot{I})$ (with $\dot{I} = I/\cK$). But the long exact sequence for the extension $0\to \cK \to I \to \dot{I} \to 0$ shows that the map $K_0(I) \to K_0(\dot{I})$ is an isomorphism (using the fact that the induced map $K_0(\cK) \to K_0(I)$ is 0).

(4) As above, we write $\dot{I},\dot{I}'$ for the quotients of $I,I'$ by $\cK$. We have \bgl\label{iso}\dot{I}=\dot{I}_1 \oplus \dot{I}_2\cong \cK\otimes C(\Tz)\oplus \cK\otimes C(\Tz)\egl and similarly for $\dot{I}'$. In particular,
$K_1(\dot{I})=K_0(\dot{I})=\Zz^2$ and the isomorphism \eqref{iso} shows that the map $\Zz^2\to \Zz^2$ induced by $\kappa$ acting on $K_1$ is the same as the map $\Zz^2\to \Zz^2$ induced by $\kappa$ on $K_0$. This latter map $\mu$ fits into the following commutative diagram of boundary maps
\bgl\label{mu}\xymatrix{K_1(C\Tz^2)=\Zz^2\ar[d]^{M^\perp\cdot}\ar[r]^\id&K_0(\dot{I}')=\Zz^2\ar[d]^\mu\\
K_1(C\Tz^2)=\Zz^2\ar[r]^{M\cdot }&K_0(\dot{I})=\Zz^2}\egl
Here, as in point (3) we identify $K_0(\dot{I}'_2),K_0(\dot{I}_2)$ with the first component and $K_0(\dot{I}'_1),K_0(\dot{I}_1)$ with the second component of $\Zz^2$. By commutativity, $\mu$ is the same as multiplication by $MM^\perp=\det M \textbf{1}$. Finally, the maps $K_1(I)\to K_1(\dot{I})$ and $K_1(I')\to K_1(\dot{I'})$ are injections so that the map $K_1I'\cong \Zz \to K_1(I)\cong \Zz$ induced by $\kappa$ is the restriction of $\mu$ to the images of these maps.
\eproof
\bremark Since $S-S=\Zz^2$, lemma \ref{Lq} shows that $\Zz^2/M\Zz^2$ is isomorphic to the quotient $S/F$ where $F=F_1+F_2$ is the subsemigroup generated by the faces $F_1$ and $F_2$.\eremark
\btheo Let $S$ be a finitely generated subsemigroup of $\Zz^2$ as above. The $K$-theory of $C^*_\lambda S$ is determined by the formula
$$K_0(C^*_\lambda S) = S/F\oplus \Zz \qquad K_1(C^*_\lambda S) = 0$$
where $F$ is the sum of the two one-dimensional faces in $S$.
\etheo
\bproof
We use the natural map $\kappa :C^*_\lambda \Nz^2 \to C^*_\lambda S$ mapping the generators of $C^*_\lambda \Nz^2 $ to $v_1,v_2\in C^*_\lambda S$ (see Lemma \ref{L19}). We then compare the long exact sequence for the extension $0\to I\to C^*_\lambda S \to C(\Tz^2)\to 0$ with the corresponding long exact sequence for the extension $0\to I'\to C^*_\lambda \Nz^2 \to C(\Tz^2)\to 0$ and use the fact that the long exact $K$-theory sequence for the second extension is explicitly known.
Using then that $K_0(I)=K_0(I')=\Zz^2$, $K_1(I')=\Zz$ and that $K_0(C\Tz^2)=\Zz^2$, $K_1(C\Tz^2)=\Zz^2$ we obtain the following morphism of exact sequences
$$\xymatrix{\ar[r]&\Zz^2\!\ar[d]\ar[r]&\Zz\ar[d]
\ar[r]&\Zz^2\ar[r]^{\beta'}\ar[d]^{\varphi}&
\Zz\ar[r]\ar[d]^{\psi}&0\ar[r]\ar[d] &\Zz^2\!\ar[r]\ar[d]&\\
\ar[r]^{M\cdot }&\Zz^2\!\ar[r]^\gamma&K_0(C^*_\lambda\! S)\!\ar[r]&
\Zz^2\ar[r]^\beta&\Zz\ar[r]&K_1(C^*_\lambda \!S)\!\ar[r]^\alpha& \Zz^2\!\ar[r]^{M\cdot}&} $$
According to Lemma \ref{L18}(2) and (4), the map $\varphi$ maps the class $[1]$ in $K_0(C\Tz^2)$ to $[1]$ and multiplies the class $b$ of the Bott element in $K_0(C\Tz^2)$ by $\det M$, while the map $\psi$ is multiplication by $\det M$. Since $\beta'(b)=1$, it follows that also $\beta (b)=1$ and of course we have $\beta ([1])=0$. Moreover $\alpha$ is 0 since the subsequent map $M\cdot$ is injective. Thus we see that $K_1(C^*_\lambda S)=0$.

It also follows that $K_0(C^*_\lambda S)$ is an extension of Im$\gamma\cong \Zz^2/M\Zz^2$ by Ker$ \beta =\Zz$ and thus that $K_0(C^*_\lambda S) = \Zz^2/M\Zz^2\oplus \Zz$. Finally the remark above shows that $\Zz^2/M\Zz^2 \cong S/F$.
\eproof

\bremark The proof of the theorem shows that the torsion part of $K_0(C^*_\lambda S)$ is generated by the classes of the projections $e_{F_1}, e_{F_2}$. \eremark

\end{document}